\documentclass{article}
\usepackage{amsthm, amsfonts, amssymb, latexsym, amsmath}

\usepackage{xcolor}
\newcommand{\F}{\mathbb{F}}
\newcommand{\Z}{\mathbb{Z}}
\newcommand{\C}{\mathbb{C}}

\title{ Binary Sequences Derived from Differences of Consecutive Primitive Roots}

\author{Arne Winterhof$^1$ and Zibi Xiao$^2$\\
~\\
$^1$ Johann Radon Institute for\\ 
Computational and Applied Mathematics\\
Austrian Academy of Sciences\\
Altenbergerstr.\ 69, 4040 Linz, Austria\\
e-mail: arne.winterhof@oeaw.ac.at\\
$^2$ College of Science\\
Wuhan University of Science and Technology\\ Wuhan 430081, Hubei, China\\
e-mail: xiaozibi@wust.edu.cn
}

\date{}

\newtheorem{lemma}{Lemma}
\newtheorem{theorem}{Theorem}

\newtheorem{proposition}{Proposition}

\begin{document}

\maketitle

\begin{abstract} Let $1<g_1<\ldots<g_{\varphi(p-1)}<p-1$ be the ordered primitive roots modulo~$p$.
We study the pseudorandomness of the binary sequence $(s_n)$ defined by $s_n\equiv g_{n+1}+g_{n+2}\bmod 2$, $n=0,1,\ldots$.
In particular, we study the balance, linear complexity and $2$-adic complexity of $(s_n)$. 
We show that for a typical $p$ the sequence $(s_n)$ is quite unbalanced. However, there are still infinitely many 
$p$ such that $(s_n)$ is very balanced. We also prove similar results for the distribution of longer patterns. Moreover, we give general lower bounds on the linear complexity and $2$-adic complexity of~$(s_n)$ and state sufficient conditions for attaining their maximums. Hence, for carefully chosen $p$, these sequences are attractive candidates for cryptographic applications. 
\end{abstract}



\section{Introduction}

For a prime $p\ge 11$, let $g_1,\ldots,g_{\varphi(p-1)}$ with $$1<g_1<g_2<\ldots<g_{\varphi(p-1)}<p-1$$
be all the primitive roots modulo $p$ in increasing order, where $\varphi(n)$ is Euler's totient
function.
The sequence $(s_n)$ derived from the parities of differences (or sums) between consecutive primitive
roots modulo $p$ is a binary sequence of period $T=\varphi(p-1)-1$ and its first period is defined by
\begin{equation}\label{sndef}
s_n\equiv g_{n+1}+g_{n+2} \bmod 2,\quad n=0,1,\ldots,T-1.
\end{equation}

Caragiu et al.\ \cite{ctkm} 
calculated the linear complexity of this sequence for the first $1000$ primes $p$  
showing that for 
$610$ primes $p$ the sequence has maximal linear complexity which may suggest this sequence for cryptography.
This has motivated us to study theoretically properties of this sequence.

{\em Balance} and uniform {\em pattern distribution} are desirable features of a cryptographic sequence.
In Section~\ref{bal} we show that the sequence $(s_n)$ is rather unbalanced if $\frac{\varphi(p-1)}{p}$ is large.
For example, if $\frac{\varphi(p-1)}{p}$ is close to its supremum $1/2$, we have for sufficiently large $p$ essentially $2T/3$
ones and $T/3$ zeros in a period of $(s_n)$. This is the case for Fermat primes $p=2^s+1$ and safe primes, that is, $(p-1)/2$ is prime. The sequence $(s_n)$ becomes more balanced with decreasing $\frac{\varphi(p-1)}{p}$. Note that for any $\varepsilon>0$ there are infinitely many primes with $\frac{\varphi(p-1)}{p}<\varepsilon$. However, for a typical
$p$ we get unbalanced sequences.
We also study the distribution of longer patterns in $(s_n)$ in Section~\ref{patt}.
 Our results on balance and pattern distribution are based on a result of Cobeli and Zaharescu on the distribution of primitive roots \cite{cz}. Note that in the special case that $p$ is either a Fermat prime or a safe prime, that is, the primitive roots coincide with the quadratic non-residues except $-1$ for the latter, the result of Ding \cite{d} on the distribution of quadratic residues can be used to improve our error term, see \cite{wx} and the Remarks below Theorem~\ref{main1}.

The {\em linear complexity} of a sequence is the length of the shortest linear feedback shift register that 
generates the sequence. A large linear complexity is essential for cryptographic applications.
For a periodic sequence $(s_n)$ of period~$T$ we can calculate the linear complexity $L(s_n)$
by
\begin{equation}\label{LS}
L(s_n)=T-\deg(\gcd(X^T-1,S(X)),
\end{equation}
where 
$$S(X)=\sum_{n=0}^{T-1}s_nX^n,$$ see for example \cite[Lemma~8.2.1]{cdr}.

The {\em $2$-adic complexity} $C(s_n)$ of a $T$-periodic binary sequence is the length of 
the shortest feedback with carry shift register and can be calculated by 
\begin{equation}\label{Cdef}
C(s_n)=\left\lfloor\log_2\left(\frac{2^T-1}{\gcd(2^T-1,S(2))}\right)\right\rfloor,
\end{equation}
where we denote by $\log_2(x)$ the binary logarithm of $x$.

For some periods $T$ any non-constant sequence of period $T$ has a large linear complexity and a large $2$-adic complexity, respectively. 
In particular, we will see in Section~\ref{lincom} that if $T=\varphi(p-1)-1$ is a prime such that $2$ is a primitive root modulo $T$ and $p\equiv 1\bmod 4$, then the linear complexity of $(s_n)$ attains its maximum $L(s_n)=T$.
Moreover, if $2^T-1$ is a Mersenne prime, then the $2$-adic complexity of $(s_n)$ attains its maximum. 

In Section~\ref{heu} we provide some experimental data which indicates that it is not difficult to find large primes $p$ such that the sequence $(s_n)$ is balanced and has a desirable pattern distribution at least for short patterns, a large linear complexity and a large $2$-adic complexity. Hence, for carefully chosen $p$ our sequences 
are attractive candidates for cryptography.

For surveys and some recent articles on linear complexity, $2$-adic complexity and related measures of pseudorandomness see \cite{cdr,homewi,howi,hu,mewi, meniwi,mewi1,mewi2,niwi,towi,wi,wx,zibi1,zibi2,xi}.

We use the notation $f(n)=O(g(n))$ if $|f(n)|\le c|g(n)|$ for some absolute constant $c>0$ and the notation $f(n)=o(g(n))$ if $g(n)\ne 0$ for sufficiently large~$n$ and $\lim\limits_{n\rightarrow\infty}\frac{f(n)}{g(n)}=0$.

\section{Balance and Pattern Distribution}

\subsection{Balance}\label{bal}

In this section we discuss the balance of the sequence $(s_n)$ of parities of differences of primitive roots modulo $p$ defined by $(\ref{sndef})$. 

\begin{theorem}\label{main1}
Let $p$ be a prime, and let $N(1)$ and $N(0)$ denote the number of $1$s and $0$s, respectively, in a period of the sequence $(s_n)$ defined by $(\ref{sndef})$ of period $T=\varphi(p-1)-1$. Then we have
$$N(1)=\left(\frac{1}{2-\varphi(p-1)/p}+o(1)\right)T$$
and $$N(0)=\left(\frac{1-\varphi(p-1)/p}{2-\varphi(p-1)/p}+o(1)\right)T,$$
where $p\rightarrow \infty$.
\end{theorem}

Proof. 
For $i\in \F_p^*$ let $c(i)=1$ if $i$ is a primitive root modulo $p$ and $c(i)=-1$ otherwise. 
For $s\ge 1$ and $\varepsilon_1, \varepsilon_2, \cdots, \varepsilon_s\in \{-1,1\}$, set
\begin{equation}\label{Neps} M(\varepsilon_1, \cdots, \varepsilon_s)
=\left| \left\{j=1,2,\ldots,p-s: c(j+i)=\varepsilon_{i+1},~i=0,\ldots, s-1\right\}\right|.
\end{equation}
Let $z=z(\varepsilon_1,\ldots,\varepsilon_s)$ be the number of $i$ with $\varepsilon_i=1$, $i=1,\ldots,s$, and put 
$$\eta=\eta(p)=\tfrac{\varphi(p-1)}{p}.$$
From \cite[Theorem~1]{cz} we get
\begin{equation}\label{czthm1}\left|M(\varepsilon_1,\ldots,\varepsilon_s)-p\eta^z
\left(1-\eta\right)^{s-z}\right|\le 2^{s-z+1}s\sqrt{p}\log p (\tau(p-1))^s,
\end{equation}
where $\tau(p-1)$ is the number of divisors of $p-1$.

Note that $\tau(p-1)=p^{O(1/\log\log p)}$, see for example \cite[Theorem 13.12]{apst}.
Then for sufficiently small $s$ with respect to $p$, $(\ref{czthm1})$ simplifies to 
\begin{equation}\label{N} M(\varepsilon_1,\ldots,\varepsilon_s)=p\eta^z
\left(1-\eta\right)^{s-z}+O\left(p^{1/2+o(1)}\right),\quad s=o(\log\log p).
\end{equation}

For a non-negative integer $k$ put
$$N_k=M(1,\underbrace{-1,\ldots,-1}_k,1),$$
that is, $z=z(1,-1,\cdots,-1,1)=2$.
$N_k$ contributes to $N(1)$ for even $k$ and to $N(0)$ for odd $k$.
Choosing
$$m=\left\lfloor \frac{\log\log p}{\log \log\log p}\right\rfloor=o(\log\log p)$$
and recall
\begin{equation}\label{phi}\frac{p}{2}> \varphi(p-1)\gg \frac{p}{\log\log p},
\end{equation}
see for example \cite[Section 18.4]{hw},
we have by $(\ref{N})$
\begin{eqnarray*}
N(1)&\ge& \sum_{k=0}^m N_{2k}=p\eta^2\sum_{k=0}^m(1-\eta)^{2k}+O\left(p^{1/2+o(1)}\right)\\
&=&\left(\frac{1}{2-\eta}+o(1)\right)T
\end{eqnarray*}
and
\begin{eqnarray*} N(0)&\ge& \sum_{k=0}^m N_{2k+1}=p\eta^2(1-\eta)\sum_{k=0}^m(1-\eta)^{2k}+O\left(p^{1/2+o(1)}\right)\\
&=&\left(\frac{1-\eta}{2-\eta}+o(1)\right)T.
\end{eqnarray*}
The result follows from these inequalities and 
$N(0)+N(1)=T$.~\hfill $\Box$\\

Remarks. 1. Large $\varphi(p-1)$:\\
We have 
$$\varphi(p-1)\le \frac{p-1}{2}$$
which is attained for {\em Fermat primes} $p$, that is, $p$ is of the form $p=2^s+1$. 
We may call $\varphi(p-1)$ {\em large} with respect to $p$ if 
$$\varphi(p-1)=\frac{p}{2}+o(p).$$ 
{\em Safe primes} $p$, that is, $(p-1)/2$
is also a {\em (Sophie Germain) prime}, are further examples of large $\varphi(p-1)=(p-3)/2$.
For primes $p$ with large $\varphi(p-1)$, a period of the sequence $(s_n)$ consists of 
$$N(0)=(1/3+o(1))T$$ 
zeros and 
$$N(1)=(2/3+o(1))T$$ 
ones and is very unbalanced. 

Note that for a Fermat prime $p$ the primitive roots modulo $p$ are exactly the quadratic non-residues and 
the proof of \cite[Theorem~3.1]{wx} can be easily modified to get the Theorem with a more precise error term.
The same applies to a safe prime $p$ for which the primitive roots modulo $p$ are the quadratic non-residues 
$\neq p-1$.

 Since the sequence $(s_n)$ is not balanced for large $\varphi(p-1)$, as in \cite{wx} we may consider the essentially  balanced sequence $(t_n)$ with $t_n=1$ whenever
$g_{n+1}=g_n+1$ and $t_n=0$ otherwise instead of $(s_n)$. \\

2. Small $\varphi(p-1)$:\\
We have $\varphi(n)\gg n/\log\log n$ which is attained for infinitely many $n$, see for example \cite{hw}.
We call $\varphi(p-1)$ of order of magnitude $p/\log \log p$ or more general with 
$$\varphi(p-1)=o(p)$$ 
{\em small}. In this case, for sufficiently large~$p$, the sequence~$(s_n)$ is essentially balanced, that is, $$N(a)=\left(\frac{1}{2}+o(1)\right)T,\quad a=0,1.$$

3. Typical $\varphi(p-1)$:\\
For an even $n$ the expected value of $\varphi(n)$ is  $4n/\pi^2$. 
More precisely, the probability that a randomly chosen even number $n$ and a random number $k$ are both divisible by a prime $r>2$
is $1/r^2$. Hence, the probability that $n$ and $k$ are co-prime is 
$$\frac{1}{2}\prod_{r>2} \left(1-\frac{1}{r^2}\right)=\frac{2}{3} \prod_r \left(1-\frac{1}{r^2}\right)=\frac{2}{3\zeta(2)} =\frac{4}{\pi^2}.$$
Here, see for example \cite[Chapters 11 and 12]{apst}, 
$$\zeta(s)=\sum_{n=1}^\infty n^{-s},\quad s\in \C,$$ 
denotes the {\em Riemann zeta function}, by Euler's product formula \cite[Theorem~11.7]{apst} we have
$$\frac{1}{\zeta(s)}=\prod_r \left(1-\frac{1}{r^s}\right),\quad {\rm Re}(s)>1,$$
and it is well-known, \cite[Theorem~ 12.17]{apst}, that $$\zeta(2)=\frac{\pi^2}{6}.$$
We call $\varphi(p-1)$ {\em typical} if 
$$\varphi(p-1)=\left(\frac{4}{\pi^2}+o(1)\right) p.$$ 
In this case we have 
$$N(1)=\left(\frac{1}{2-4/\pi^2}+o(1)\right)T=(0.627\ldots+o(1))T$$ and
$$N(0)=\left(\frac{1-4/\pi^2}{2-4/\pi^2}+o(1)\right)T=(0.372\ldots+o(1))T.$$

\subsection{Pattern Distribution}\label{patt}

Now we extend Theorem~\ref{main1} to longer patterns of fixed length $\ell$ and $p\rightarrow \infty$.
\begin{theorem}
Let $(a_0,\ldots,a_{\ell-1})\in \F_2^\ell$ be a pattern of fixed length $\ell\ge 1$ with $w$ coordinates equal to $1$ and $\ell-w$ coordinates equal to $0$. Let $N_\ell(w)$ be the number of $n=0,1,\ldots,T-\ell$
with $s_{n+i}=a_i$ for $i=0,\ldots,\ell-1$.
Then we have
$$N_\ell(w)=\left(\left(\frac{1}{2-\varphi(p-1)/p}\right)^w\left(\frac{1-\varphi(p-1)/p}{2-\varphi(p-1)/p}\right)^{\ell-w}+o(1)\right)T,\quad p\rightarrow \infty.$$
\end{theorem}

Proof. Without loss of generality we consider 
$(a_0,\ldots,a_{\ell-1})=(\underbrace{1,\ldots,1}_w,\underbrace{0,\ldots,0}_{\ell-w})$.

Recall $(\ref{Neps})$ and put 
\begin{eqnarray*}&&N_{k_1,\ldots,k_\ell}=\\
&&M(1,\underbrace{-1,\ldots,-1}_{2k_1},1,\ldots,1,\underbrace{-1,\ldots,-1}_{2k_w},1,\underbrace{-1,\ldots,-1}_{2k_{w+1}+1},1,\ldots,
1,\underbrace{-1,\ldots,-1}_{2k_\ell+1},1).
\end{eqnarray*}
Put 
$$m=\left\lfloor \frac{\log\log p}{\log\log\log p}\right\rfloor.$$
Then we have
\begin{eqnarray*} N_\ell(w)&\ge& \sum_{k_1,\ldots,k_\ell=0}^m N_{k_1,\ldots,k_\ell}\\
&\ge& p\eta^{\ell+1}(1-\eta)^{\ell-w}\sum_{k_1,\ldots,k_\ell=0}^m(1-\eta)^{2(k_1+\ldots+k_\ell)}+O\left(p^{1/2+o(1)}\right)
\end{eqnarray*}
by $(\ref{N})$, where $\eta=\frac{\varphi(p-1)}{p}$
and thus
\begin{eqnarray*}
N_\ell(w)
 &\ge& T \eta^\ell (1-\eta)^{\ell-w} \left(\sum_{k=0}^m (1-\eta)^{2k}\right)^\ell+O\left(p^{1/2+o(1)}\right)\\
 &=& T \eta^\ell (1-\eta)^{\ell-w} \left(\frac{1-(1-\eta)^{2(m+1)}}{1-(1-\eta)^2}\right)^\ell+O\left(p^{1/2+o(1)}\right)\\
&=& T(1-\eta)^{\ell-w}\left(\frac{1+o(1)}{2-\eta}\right)^\ell+O\left(p^{1/2+o(1)}\right)\\
&=&T\left((1-\eta)^{\ell-w}\left(\frac{1}{2-\eta}\right)^\ell+o(1)\right).
\end{eqnarray*}
 In the last step we used $T\gg p/\log\log p$ and $0<\frac{1}{2-\eta}<1$ since $0<\eta<1/2$, see the remark after Theorem~\ref{main1}. We recall that $\ell$ is fixed.

We have ${\ell \choose w}$ patterns 
$(a_0,\ldots,a_{\ell-1})$ with $w$ coordinates equal 
to $1$.
Since 
$$\sum_{w=0}^\ell {\ell \choose w}\left(\frac{1}{2-\eta}\right)^w\left(\frac{1-\eta}{2-\eta}\right)^{\ell-w}=
\left(\frac{1}{2-\eta}+\frac{1-\eta}{2-\eta}\right)^\ell=1$$
the main term of this lower bound is optimal and the result follows.
~\hfill $\Box$\\

Remark. For large $\varphi(p-1)=\frac{p}{2}+o(p)$ we get
$$N_\ell(w)=\left(\left(\frac{2}{3}\right)^w\left(\frac{1}{3}\right)^{\ell-w}+o(1)\right)T.$$
For small $\varphi(p-1)=o(p)$ we get
$$N_\ell(w)=\left(\left(\frac{1}{2}\right)^\ell+o(1)\right)T$$
and for typical $\varphi(p-1)=\frac{4p}{\pi^2}+o(p)$ we have
$$N_\ell(w)=\left((0.627\ldots)^w (0.372\ldots)^{\ell-w}+o(1)\right)T.$$

\section{Linear Complexity and $2$-Adic Complexity}
\subsection{Linear Complexity} 
\label{lincom}

In this section, we estimate the linear complexity of the $T$-periodic sequence~$(s_n)$ defined by $(\ref{sndef})$.
In particular, we give a sufficient condition for attaining the maximal value $L(s_n)=T$.

For integers $m$ and $q$ with $\gcd(m,q)=1$ we denote by ${\rm ord}_m(q)$ the {\em order} of~$q$ modulo $m$.
Note that $\varphi(p-1)$ is even for $p\ge 5$, that is, $T=\varphi(p-1)-1$ is odd and $T\ge 3$ for $p\ge 11$.
\begin{proposition}\label{main2}
Let $p$ be a sufficiently large prime and $T=\varphi(p-1)-1$. Let $T=p_1^{e_1}\cdots p_r^{e_r}$ be the prime factorization of $T$ with pairwise distinct odd primes $p_1,\ldots,p_r$ and $e_i\ge 1$ for $i=1,\ldots,r$.
Then the linear complexity of the sequence~$(s_n)$ of period $T$ defined by~$(\ref{sndef})$ satisfies
$$L(s_n)\ge \min\left\{{\rm ord}_{p_1}(2),\ldots,{\rm ord}_{p_r}(2)\right\}+\varepsilon,$$
where
\begin{equation}\label{eps} \varepsilon=\left\{\begin{array}{cc} 1, & p\equiv 1\bmod 4,\\ 0, & p\equiv 3\bmod 4.\end{array}\right.
\end{equation}
In particular, if $T$ is a prime and $2$ is a primitive root modulo~$T$, then 
$$L(s_n)\left\{\begin{array}{lc} =T,& p\equiv 1\bmod 4,\\
\ge T-1, & p\equiv 3\bmod 4.\end{array}\right.$$ 
\end{proposition}

The proof is based on a slightly more precise version of \cite[Theorem 3.3.1]{cdr}.

\begin{lemma}\label{ding}
 Let $T=p_1^{e_1}\cdots p_r^{e_r}$ be the prime factorization of an odd integer $T\ge 3$ with pairwise distinct primes $p_1,\ldots,p_r$ and $e_i\ge 1$ for $i=1,\ldots,r.$
 Then for each non-constant sequence $(s_n)$ over $\F_2$ of period $T$ we have
 $$L(s_n)\ge \min\left\{{\rm ord}_{p_1}(2),\ldots,{\rm ord}_{p_r}(2)\right\}+S(1),$$
 where 
 $$S(1) = \sum_{n=0}^{T-1}s_n\in \F_2=\{0,1\}.$$
\end{lemma}
Proof.  
Since $T$ is odd we have $\gcd(X-1,X^{T-1}+\ldots+X+1)=1$
and thus 
$$\gcd(X^T-1,S(X))=\gcd(X-1,S(X))\gcd(X^{T-1}+\ldots+X+1,S(X)).$$ 
From the proof of \cite[Theorem 3.3.1]{cdr} we know that 
$$T-\deg(\gcd(X^{T-1}+\ldots+X+1,S(X)))\ge \min\{{\rm ord}_{p_1}(2),\ldots,{\rm ord}_{p_r}(2)\}.$$
Now $\gcd(X-1,S(X))=X-1$ if $S(1)=0$ and $\gcd(X-1,S(X))=1$ if $S(1)=1$ and the result follows from $(\ref{LS})$.~\hfill $\Box$\\

Now we study the value of $S(1)$.

\begin{lemma}
 For a prime $p\equiv 1\bmod 4$ and the sequence $(s_n)$ defined by $(\ref{sndef})$ we have 
 $$S(1)=1.$$
\end{lemma}

Proof. By the definition of $(s_n)$ we have
$$S(1)=\sum_{n=0}^{\varphi(p-1)-2} s_n= \sum_{n=0}^{\varphi(p-1)-2} (g_{n+1}+g_{n+2})= g_1+g_{\varphi(p-1)} \in \F_2.$$
For an arbitrary primitive root $g$ modulo $p$ we have $g^{(p+1)/2}\equiv -g\bmod p$.
Since $\gcd((p+1)/2,p-1)=1$ for $p\equiv 1\bmod 4$, it follows that $-g$ is also
a primitive root modulo~$p$. This shows that if $g_1$ denotes the smallest primitive root modulo~$p$,
then $p-g_1$ is the largest primitive root modulo $p$, that is, $g_1+g_{\varphi(p-1)}=p$ in~$\Z$.
Thus we have $$g_1+g_{\varphi(p-1)}=1 \in \F_2,$$
which completes the proof.~\hfill $\Box$\\

Remark. 
For $p\equiv 3\bmod 4$ both possible values of $S(1)$ can be attained.
For example, $S(1)= 2+8= 0\in \F_2$ for $p=11$ and $S(1)= 2+15= 1 \in \F_2$ for $p=19$.\\

For proving Proposition~\ref{main2} it remains to verify that the sequence $(s_n)$ defined by $(\ref{sndef})$ is non-constant for a sufficiently large prime $p$. By $(\ref{phi})$ and Theorem~\ref{main1} we have
$$N(1)\ge \left(\frac{1}{2}+o(1)\right)T \quad \mbox{and}\quad N(0)\ge \left(\frac{1}{3}+o(1)\right)T.$$
Hence, $N(1)$ and $N(0)$ are both positive for sufficiently large $p$ and $(s_n)$ is not constant. Hence,  Lemma~\ref{ding} is applicable and completes the proof of Proposition~\ref{main2}.~\hfill $\Box$

\subsection{$2$-Adic Complexity}

Now we estimate the $2$-adic complexity of $(s_n)$ defined by $(\ref{sndef})$.
\begin{proposition}
Let $p$ be a sufficiently large prime and $T=\varphi(p-1)-1$. Let~$q$ be the smallest prime divisor of $2^T-1$.
Then the $2$-adic complexity of the sequence~$(s_n)$ of period $T$ defined by~$(\ref{sndef})$ satisfies
$$C(s_n)\ge \lfloor\log_2(q)\rfloor.$$
In particular, if $2^T-1$ is a (Mersenne) prime, then $$C(s_n)=\lfloor\log_2(2^T-1)\rfloor.$$ 
\end{proposition}

Since $(s_n)$ is not constant for sufficiently large $p$, by Theorem~\ref{main1} it is enough to
verify the following lemma, which may be of independent interest.

\begin{lemma}
 Let $q$ be the smallest prime divisor of $2^T-1$.
 Then for each non-constant sequence $(s_n)$ over $\F_2$ of period $T$ we have
 $$C(s_n)\ge \lfloor \log_2(q)\rfloor.$$
\end{lemma}
Proof. Put $d=\gcd(S(2),2^T-1)$. We have $d=2^T-1$ if and only if $S(2)\in \{0,2^T-1\}$, that is, $(s_n)$ is constant.

Now assume that $(s_n)$ is not constant and $q$ denotes the smallest prime divisor of $2^T-1$.
Then we have $d\le \frac{2^T-1}{q}$ and thus
$$C(s_n)=\left\lfloor\log_2\left(\frac{2^T-1}{\gcd(S(2),2^T-1)}\right)\right\rfloor\ge \lfloor\log_2(q)\rfloor$$
by $(\ref{Cdef})$.~\hfill $\Box$\\

Remark. Note that there are highly predictable sequences with both maximum linear complexity and maximum $2$-adic complexity, for example, 
any sequence with only one non-zero entry in a period. Hence, studying the balance and pattern distribution is always a must to test a sequence for 
suitability in cryptography.

\section{Heuristic}\label{heu}

To guarantee a rather balanced sequence with large linear complexity and large $2$-adic complexity we need primes $p$ such that
\begin{itemize}
    \item The ratio $\frac{\varphi(p-1)}{p}$ is small.
    \item The period $T=\varphi(p-1)-1$ contains only large prime divisors $q$ such that ${\rm ord}_q(2)$ is also large. This is guaranteed if $T$ is prime and $2$ is a primitive root modulo $2$.
    \item The Mersenne number $2^T-1$ contains only large prime divisors. This is guaranteed if $2^T-1$ is a Mersenne prime.
\end{itemize}

In the following table we list primes $T$ for which $2^T-1$ is a Mersenne prime and the largest primes $p$
with $T=\varphi(p-1)-1$. For these primes we have $L(s_n)\ge {\rm ord}_T(2)+S(1)$ with $S(1)$ defined by $(\ref{eps})$, and $C(s_n)$ is maximal.\\

    \begin{tabular}{c|c|c|c}
$T$ & $p$ & ${\rm ord}_T(2)$ &  $\frac{\varphi(p-1)}{p}$\\ \hline
       $3$  & $13$ & $2$ & $\frac{4}{13}=0.307\ldots$ \\\
       $5$ & $19$ & $4$ &$\frac{6}{19}=0.315\ldots$ \\
       $7$ & $31$ & $3$ &$\frac{8}{31}=0.258\ldots$ \\
       $19$ & $67$& $18$ &$\frac{20}{67}=0.298\ldots$ \\
       $31$ & $103$ & ${\bf 5}$ &$\frac{32}{103}=0.310\ldots$\\
       $107$ & $379$ & $106$ &$\frac{108}{379}=0.284\ldots$ \\
       $127$ & $409$ & ${\bf 7}$ &$\frac{128}{409}=0.312\ldots$ \\
       $1279$ & $5281$& $639$ &$\frac{1280}{5281}=0.242\ldots$\\
       $2203$ & $6619$ & $734$ & $\frac{2204}{6619}=0.331\ldots$
    \end{tabular}\\

Now we also list some primes $T$ for which $2^T-1$ is not a prime. We denote by~$q$ the smallest prime divisor of $2^T-1$ from which we can derive the lower bound $C(s_n)\ge \lfloor \log_2(q)\rfloor$ on the $2$-adic complexity.\\

   \begin{tabular}{c|c|c|c|c|c}
$T$ &  $q$ & $\lfloor\log_2(q)\rfloor$& $p$ &${\rm ord}_T(2)$ &  $\frac{\varphi(p-1)}{p}$\\ \hline
       $11$ & $23$     & $4$ & $43$  & $10$ & $0.279\ldots$ \\
       $23$ & $47$     & $5$ & $79$ & $11$ &$0.303\ldots$ \\
       $43$ & $431$    & $8$ & $139$         & $14$ &$0.316\ldots$ \\
       $47$ & $2351$   & $11$ & $211$ & $23$ &$0.227\ldots$ \\
       $53$ & $6361$   & $12$ & $163$         & $52$ &$0.331\ldots$\\
       $59$ & $179951$ & $17$ & $199$         & $58$ & $0.301\ldots$\\
       $71$ & $228479$ & $17$ & $271$ & $35$ &$0.265\ldots$ \\
       $79$ & $2687$   & $11$ & $331$         & $39$ &$0.209\ldots$ \\
       $83$ & $167$    & ${\bf 7}$ & $197$     & $82$ & ${\bf 0.426}\ldots$\\
      $131$ & $263$    & ${\bf 8}$ & $269$         & $130$ & ${\bf 0.490}\ldots$ \\
      $163$ & $150287$ & $17$ & $499$ & $162$ & $0.328\ldots$ \\
      $167$ & $2349023$ & $21$ & $523$ & $83$ & $0.321\ldots$\\
      $179$ & $359$ & ${\bf 8}$& $419$ & $178$ & ${\bf 0.429\ldots}$\\
      $191$ & $383$ & ${\bf 8}$ &  $673$ & $95$ & $0.285\ldots$\\
      $199$ & $164504919713$ & $37$ & $751$ & $99$ & $0.255\ldots$
    \end{tabular}\\

We may consider the following features undesirable and emphasized this in the tables (boldface):
\begin{itemize}
\item The value $\lfloor\log_2(q)\rfloor$ is small, say, smaller than $\frac{T}{10}$. Then a very large $2$-adic complexity cannot be guaranteed.
\item The order of $2$ modulo $T$ is small, say, smaller than $\frac{T}{4}$. Then a very large linear complexity cannot be guaranteed.
\item The ratio $\frac{\varphi(p-1)}{p}$ is large, say, at least $\frac{1}{3}$. Then for sufficiently large $p$ the sequence contains at least $60$ percent ones and is rather unbalanced. Moreover, the frequency of the pair $11$ is at least $36$ percent whereas the frequency of $00$ is at most $16$ percent of the period.     
\end{itemize}    
Still it seems to be not difficult to find large primes $T$ and $p$ with $T=\varphi(p-1)-1$ without these undesirable features.

\section*{Acknowledgment}
The first author was partially supported by the Austrian Science Fund FWF Project P 30405-N32. 
The second author was supported in part by the Chinese Scholarship Council and in part by the National Natural Science Foundation of China under Grant 12061027.

\end{document}